\documentclass[letterpaper, 10 pt, conference]{ieeeconf}  

\IEEEoverridecommandlockouts                              

\usepackage{graphics} 
\usepackage{epsfig} 
\usepackage{mathptmx} 
\usepackage{times} 
\usepackage{amsmath} 
\usepackage{amssymb}  

\usepackage{color}
\newtheorem{remark}{Remark}
\newtheorem{theorem}{Theorem}
\newtheorem{lemma}{Lemma}
\newtheorem{proposition}{Proposition}

\DeclareMathAlphabet{\cmcal}{OMS}{cmsy}{m}{n}

\title{\LARGE \bf
Positive Regulation of an Axial Dispersion Tubular Reactor by Optimal Boundary Feedback Control*
}

\author{Violaine Piengeon$^{1}$ and Joseph J. Winkin$^{1}$
\thanks{*Violaine Piengeon is a Research Fellow of the Fonds de la Recherche Scientifique - FNRS.}
\thanks{$^{1}$University of Namur, Department of Mathematics \&  Namur Research Institute
for Complex Systems (naXys), Rue Graf\'e 2, B-5000 Namur, Belgium
            {\tt\small violaine.piengeon@unamur.be, joseph.winkin@unamur.be}}%
}

\begin{document}

\maketitle
\thispagestyle{empty}
\pagestyle{empty}

\begin{abstract}
The dynamics of the reactant in an isothermal axial dispersion tubular reactor is considered. Our goal consists in designing the inlet reactant concentration, under the form of a state feedback law, in order to modify the state behavior according to relevant process requirements, while keeping the positivity and stability properties of the uncontrolled dynamics. The designed control is also constrained to be positive to respect physical requirements. First, we state and analytically solve a parameterized optimal feedback design problem. 
Next, we address a positive observer-based feedback stabilization problem.
The results are illustrated by means of numerical simulations.
\end{abstract}

\section{INTRODUCTION}
A key question in process control is how to regulate a system to fit performance requirements whilst respecting physical constraints inherent to the process. For chemical concentration dynamics, the  nonnegativity of state trajectories is clearly a fundamental requirement.
This motivates the use of positive distributed parameter systems theory (e.g. \cite{Laabissi2001,Laabissi2003}), which provides a mathematical framework for analyzing systems whose state trajectories remain in the positive cone of an appropriate Hilbert space whenever the initial conditions are nonnegative. Unfortunately, using classical exact or observer-based state feedback techniques (see e.g \cite{Curtain_Zwart2020}) on positive open-loop systems does not guarantee to respect the nonnegativity requirement. The latter is therefore an essential additional constraint to take into account in the design of feedback controllers, especially for chemical processes. 

When chemical processes occur in tubular reactors, that cover a large class of applications, the dynamics is often modeled by reaction-convection-diffusion (RCD) systems; see e.g. \cite{Laabissi2001,Laabissi2003,Dochain2025,Delattre2004,Winkin2000}. Feedback design problems taking into account positivity have been studied for infinite-dimensional systems with bounded control operators notably in \cite{Abouzaid2010,Achhab2017,Abouzaid2022}, and for RCD systems in \cite{Yevgenieva2026}. However, in the references mentioned above, there is neither specific optimal control design, nor observer-based feedback design. Overall, to the best of our knowledge, these design problems have been rarely considered for positive infinite-dimensional systems. The recent contribution \cite{PiengeonWinkin2026} can be seen as an attempt in that direction, where such problems are studied for a one-dimensional pure diffusion system.

The main goal of this paper is to design positivity preserving feedback controllers for the reactant concentration dynamics in an isothermal axial dispersion tubular reactor \cite[Subsection 1.4.1]{Dochain2025}. First, a relevant parametrized optimal feedback design problem, inspired by \cite[Section 3]{PiengeonWinkin2026}, is formulated, well-motivated and analytically solved. The motivations include mathematical and application-oriented interpretations and the derivation of the solution in closed-form relies on a complete spectral analysis of the dynamics generator. Then, a positive observer-based feedback stabilization problem is addressed. More precisely, we report an effective characterization of observer-based boundary feedback laws that ensure stabilization of both the controlled system and the proposed state observer with boundary output error injection. An educated guess for an effective characterization of positivity is also obtained by analyzing a spatial discretization of an equivalent (in terms of positivity) simpler system.

The paper is organized as follows. In Section~\ref{section:model_problem}, the considered model is introduced and the challenges are explained. Section~\ref{section:exact_feedback} is dedicated to the parameterized optimal feedback design problem, while Section~\ref{section:observer_based_feedback} addresses the observer-based boundary feedback design problem. The educated guess for positivity and numerical simulations is provided in Section~\ref{section:numerical_simulations}, followed by concluding remarks in Section~\ref{section:conclusion}.

The following notations will be used throughout. $\mathbb{N}$ and $\mathbb{R}$ denote the sets of nonnegative integers and real numbers, respectively. The set $\mathrm{H}^2(0,1)$ denotes the Sobolev space of functions belonging to the Lebesgue space $\mathrm{L}^2(0,1)$ whose first and second weak derivatives are also in $\mathrm{L}^2(0,1)$. For any linear operator $A$ whose spectrum reduces to the point spectrum,
its spectral bound is given by $s(A) := \sup \; \textbf{\rm Re} \; \sigma_p(A)$, where the set $\sigma_p(A)$ denotes its point spectrum.

\section{MODEL AND PROBLEM FORMULATION}\label{section:model_problem}
Let us consider an isothermal axial dispersion tubular reactor with linear kinetics. 
This process kinetics can be characterized by the chemical reaction $C\longrightarrow bP$, where $C$ is the reactant, $P$ is the product and $b$ is the stoichiometric coefficient of the reaction; see e.g. \cite[Subsection 1.4.1]{Dochain2025}.
Moreover, the dynamics of the reactant 
is described by the following reaction-convection-diffusion dimensionless model \cite[Subsection 1.4.1]{Dochain2025}, \cite[Section 2]{Delattre2004}: 
\begin{equation}\label{eq:system}
    \frac{\partial c}{\partial \tau}(\zeta,\tau) = \frac{1}{Pe}\frac{\partial^2 c}{\partial {\zeta}^2}(\zeta,\tau)-\frac{\partial c}{\partial {\zeta}}(\zeta,\tau)-rc(\zeta,\tau),
\end{equation}
with Danckwerts’ boundary conditions 
\begin{equation}\label{eq:boundary_cdts}
    \dfrac{1}{Pe}\dfrac{\partial c}{\partial \zeta}(0,\tau)-c(0,\tau) = -c_{in}(\tau),\quad
    \dfrac{1}{Pe}\dfrac{\partial c}{\partial \zeta}(1,\tau) = 0,
\end{equation}
where the dimensionless variables $\zeta\in[0,1]$ and $\tau\geq0$, and the (dimensionless) Peclet $Pe$ and Damköhler $r$ numbers are defined by 
\begin{equation*}
    \zeta = \frac{z}{L},\quad \tau = \frac{vt}{L},\quad Pe = \frac{vL}{D_a}, \quad r=\frac{k_0L}{v},
\end{equation*}
respectively, where $z\in[0,L]$ is the spatial variable, $t\geq0$ is the time, and $D_a$, $v$, $k_0$, and $L$ are positive parameters denoting the axial dispersion coefficient ($m^2/s$), the fluid superficial velocity ($m/s$), the kinetic constant ($s^{-1}$), and the reactor length ($m$), respectively. Here, $c(\zeta,\tau)$ is the concentration of the reactant ($mol/l$) and $c_{in}(\tau)$ denotes the inlet reactant concentration ($mol/l$).
Within our theoretical control framework, $c(\tau):=c(\cdot,\tau)$ is the state (profile) at time $\tau$ and $c_{in}$ denotes the boundary control. The concentration of the reactant is measured at the outlet of the tubular reactor, so that the measured output $y$ that is considered consists in a boundary observation given by
\begin{equation}\label{eq:output}
    y(\tau) = c(1,\tau), \quad \tau\geq 0.
\end{equation}

Without any control, that is when $c_{in}=0$, 
it is well-known that the system \eqref{eq:system}-\eqref{eq:boundary_cdts} is exponentially stable and positive as the associated operator $A_{0}:D(A_0)\subset\mathrm{L}^2(0,1)\rightarrow \mathrm{L}^2(0,1)$, defined by $A_{0} = \frac{1}{Pe}\frac{{\rm  d}^2}{{\rm  d}  \zeta^2}-\frac{{\rm  d}}{{\rm  d}  \zeta}-rI$ with domain
\begin{equation*}
    \begin{aligned}
    D(A_0) = \biggl\{c\in \mathrm{H}^2(0,1) : &\frac{1}{Pe}\frac{{\rm  d}  c}{{\rm  d}  \zeta}(0)-c(0) =  0, \\
    &\frac{1}{Pe}\frac{{\rm  d}  c}{{\rm  d}  \zeta}(1) =0 \biggr\},
    \end{aligned}
\end{equation*}
generates an exponentially stable and positive $C_0$-semigroup on the state space $\mathrm{L}^2(0,1)$; see \cite{Winkin2000}, \cite{Laabissi2001}. The stability property means that, without any action on the tubular reactor, the reactant will progressively disappear along time, at a fixed rate depending exclusively upon the system parameters, by being transformed into product. Positivity ensures that the state remains nonnegative, in agreement with the physical nature of concentration as a nonnegative quantity.

In the following two sections, the goal is to design the inlet reactant concentration $c_{in}$, under the form of a feedback law, in order to modify the reactant behavior according to relevant process requirements, while keeping the positivity and stability properties. The designed control is also constrained to be nonnegative as it represents a concentration, that is $c_{in}$ must satisfy $c_{in}(\tau)\geq0$ for all $\tau\geq0$. 

Suppose now that $c_{in}$ is an exact state feedback law of the form 
\begin{equation}\label{eq:feedback_exact}
    c_{in}(\tau) = \kappa c(0,\tau),\quad \tau\geq0,
\end{equation}
where $\kappa\in\mathbb{R}$ denotes the feedback gain, to be designed. With \eqref{eq:feedback_exact}, the control strategy consists in adding along time an amount of reactant which is proportional to the concentration of the reactant already present at the inlet of the tubular reactor. Plugging \eqref{eq:feedback_exact} into \eqref{eq:system}-\eqref{eq:boundary_cdts} gives a closed-loop system that may be written in the following abstract form  
\begin{equation}\label{eq:closed_loop}
    \dot{c}(\tau) = A_{\kappa}c(\tau),
\end{equation}
where $A_{\kappa}:D(A_{\kappa})\subset\mathrm{L}^2(0,1)\rightarrow \mathrm{L}^2(0,1)$ is the operator defined by
\begin{equation}\label{eq:A_kappa}
    A_{\kappa} = \frac{1}{Pe}\frac{{\rm  d}^2}{{\rm  d}  \zeta^2}-\frac{{\rm  d}}{{\rm  d}  \zeta}-rI
\end{equation}
with domain
\begin{equation}\label{eq:domain_A_kappa}
    \begin{aligned}
    D(A_{\kappa}) = \biggl\{c\in \mathrm{H}^2(0,1) : &\frac{1}{Pe}\frac{{\rm  d}  c}{{\rm  d}  \zeta}(0)+(\kappa-1)c(0) =  0, \\
    &\frac{1}{Pe}\frac{{\rm  d}  c}{{\rm  d}  \zeta}(1) =0 \biggr\}.
    \end{aligned}
\end{equation}
Since $-A_{\kappa}$ with domain $D(-A_{\kappa}):=D(A_{\kappa}),$ is a Sturm-Liouville operator, $A_{\kappa}$ generates a $C_0$-semigroup $(T_{A_{\kappa}}(t))_{t\geq0}$ on $\mathrm{L}^2(0,1)$.

A classical design strategy in control theory therefore involves selecting a low (feedback) gain to avoid a huge perturbation of the dynamics that may have undesirable consequences like an actuator saturation; see \cite{Lin2009LowGA}. In our framework, this low gain type problem is constrained and may be formalized in the following form:
\begin{subequations}\label{eq:pbm_opti_without_alpha}
    \begin{align}
        \displaystyle \min_{\kappa \,\in \,\mathbb{R}}\quad &\kappa^2, \label{eq:pbm_opti_without_alpha_obj}\\
        \text{s.t.}
        \quad & (T_{A_{\kappa}}(t))_{t\geq0} \text{ is }\text{exponentially stable},\label{eq:pbm_opti_without_alpha_stability}\\
        \quad & (T_{A_{\kappa}}(t))_{t\geq0} \text{ is positive,} \label{eq:pbm_opti_without_alpha_system_positivity}
        \\
        \quad & \kappa\geq0. \label{eq:pbm_opti_without_alpha_input_positivity}
    \end{align}
\end{subequations}
In the sequel, the optimal control problem \eqref{eq:pbm_opti_without_alpha} will be denoted by $\mathbf{P}$. Observe that condition \eqref{eq:pbm_opti_without_alpha_input_positivity} when combined with \eqref{eq:pbm_opti_without_alpha_system_positivity} allows to ensure the required nonnegativity property for $c_{in}$. Obviously, since conditions \eqref{eq:pbm_opti_without_alpha_stability} to \eqref{eq:pbm_opti_without_alpha_input_positivity} are satisfied for $\kappa=0$, the unique solution $\kappa^*$ of the $\mathbf{P}$ problem \eqref{eq:pbm_opti_without_alpha} is trivial as $\kappa^*=0$, i.e. the optimal control design strategy consists in doing nothing. That is the reason why in the following section, a positive parameter $\alpha$ is relevantly introduced in the control design problem. 

\section{OPTIMAL FEEDBACK}\label{section:exact_feedback}
Let us denote the growth bound of $(T_{A_{\kappa}}(t))_{t\geq0}$ by $\omega_0(\kappa)$. Condition \eqref{eq:pbm_opti_without_alpha_stability} is equivalent to $\omega_0(\kappa)<0$; see \cite{Curtain_Zwart2020}. For any $\alpha>0$, the following optimal control problem $\mathbf{P_{\alpha}}$ is now considered for the system \eqref{eq:system}-\eqref{eq:boundary_cdts} under the feedback control law \eqref{eq:feedback_exact}: 
\begin{subequations}\label{eq:pbm_opti}
    \begin{align}
        \displaystyle \min_{\kappa \,\in \,\mathbb{R}}\quad &\kappa^2, \label{eq:pbm_opti_obj}\\
        \text{s.t.}
        \quad & -\alpha\leq\omega_0(\kappa)<0,\label{eq:pbm_opti_stability}\\
        \quad & (T_{A_{\kappa}}(t))_{t\geq0} \text{ is positive,} \label{eq:pbm_opti_system_positivity}
        \\
        \quad & \kappa\geq0. \label{eq:pbm_opti_input_positivity}
    \end{align}
\end{subequations}
\begin{remark}
    (i) In comparison with the $\mathbf{P}$ problem, the $\mathbf{P_{\alpha}}$ problem is parametrized by $\alpha>0$ via the left-hand side of \eqref{eq:pbm_opti_stability}. The additional inequality $-\alpha\leq\omega_0(\kappa)$ allows to constrain the value of $\omega_0(\kappa)$. Condition \eqref{eq:pbm_opti_stability} therefore means that the $C_0$-semigroup $(T_{A_{\kappa}}(t))_{t\geq0}$ is indeed exponentially stable, but with a decay rate that is not faster than $\alpha$. If $\alpha$ is chosen smaller than the opposite of the growth bound associated with the exponentially stable $C_0$-semigroup $(T_{A_{0}}(t))_{t\geq0}$, then the control strategy with the feedback law~\eqref{eq:feedback_exact} and the $\mathbf{P_{\alpha}}$ problem consists in slowing down the reactant consumption. Moreover, the smaller $\alpha$ is, the slower the reaction is intended to be. This control strategy may be useful if the process requirements consists in maintaining a minimum amount of reactant for a longer time, particularly when the reactant is required to sustain a downstream process or a prescribed production rate. \\
    (ii) Consider the closed-loop transfer function ${\hat g}_{cl}$ (of system~\eqref{eq:system}-\eqref{eq:boundary_cdts} under the feedback~\eqref{eq:feedback_exact}) from an external input $v(t)$ acting at $z=0$, such that $c_{in}(t) = \kappa c(0,t)+v(t)$, to the output \eqref{eq:output}. The static gain ${\hat g}_{cl}(0)$ is given by
    \begin{equation*}
        {\hat g}_{cl}(0) = \frac{1}{\left(p-\frac{1}{Pe}q\right)-\kappa p},
    \end{equation*}
    where 
    \begin{equation*}
        p := \frac{a_+e^{-a_-}-a_-e^{-a_+}}{a_+-a_-}, \quad q := a_+a_-\frac{e^{-a_-}-e^{-a_+}}{a_+-a_-},
    \end{equation*} 
    and $a_{\pm} := (Pe\pm\sqrt{Pe^2+4rPe})/2$. Observe that ${\hat g}_{cl}(0)$ is well defined for all $\kappa\in\mathbb{R}, \kappa\neq\kappa_0$, where
    \begin{equation}\label{eq:kappa_0}
        \kappa_0:=1-\frac{a_+a_-}{Pe}\left(\frac{e^{-a_-}-e^{-a_+}}{a_+e^{-a_-}-a_-e^{-a_+}}\right).
    \end{equation}
    Since $a_+>0$ and $a_-<0$, it follows that $p>0$, $q<0$ and $\kappa_0>1$. Consequently, ${\hat g}_{cl}(0)$, viewed as a function of $\kappa$, is increasing and positive on the interval $[0,\kappa_0)$, as illustrated in Figure \ref{fig:static_gain} for $Pe=2$ and $r=1$.
    \begin{figure}[thpb]
      \centering
      \vspace{-0.4cm}
      \includegraphics[scale=0.6]{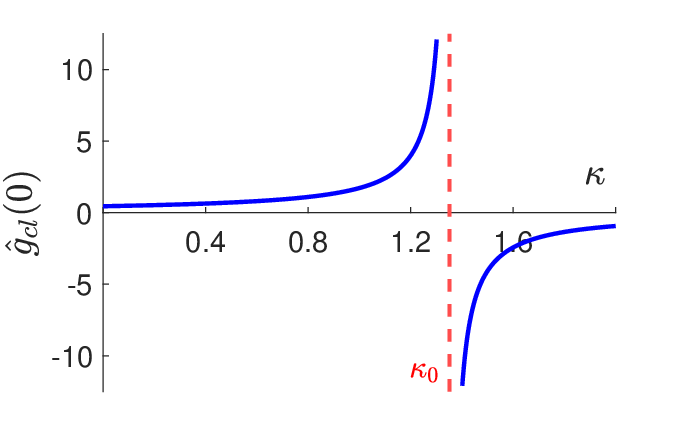}
      \vspace{-0.7cm}
      \caption{Static gain ${\hat g}_{cl}(0)$ as a function of $\kappa$ ($Pe=2, r=1$)}
      \label{fig:static_gain}
   \end{figure}
    In Theorem \ref{theorem:solution_pbm_opti}, we will prove that the feasible set of the $\mathbf{P_{\alpha}}$ problem is included in $[0,\kappa_0)$. Therefore, $\mathbf{P_{\alpha}}$ can be interpreted as minimizing the static gain, i.e. the (steady) asymptotic value of the step response (representing the concentration of the reactant at the outlet of the tubular reactor) of the closed-loop system, under the constraints \eqref{eq:pbm_opti_stability}--\eqref{eq:pbm_opti_input_positivity}. From a process control point of view, the main goal of such a design could be seen as maximizing reactant conversion and minimizing the amount of unreacted reactant leaving the reactor.
\end{remark}

In the light of those comments, the $\mathbf{P_{\alpha}}$ problem may be viewed as a trade-off between a sufficiently slow decay rate and sufficiently low feedback and static gains, while maintaining nonnegativity and exponential stability of the state trajectories and nonnegativity of the control $c_{in}$. 

With a view to solving the $\mathbf{P_{\alpha}}$ problem, for any $\alpha>0$, a complete spectral analysis of $A_{\kappa}$ is performed in Lemma \ref{lemma:spectral_analysis}. The presentation and arguments of the latter are inspired by \cite[Section 3.3]{Yevgenieva2026}. For notational simplicity, let $\bar{\kappa}$ and $u_n, n\in\mathbb{N},$ denote the following real numbers
\begin{equation*}
    \bar{\kappa}:=\frac{1}{2}+\frac{1}{Pe+2},\quad u_n := -\frac{\pi^2(1+2n)^2}{2Pe^2}+\frac{1}{2}.
\end{equation*}
Observe that $u_0<\bar{\kappa}<1<\kappa_0$, where $\kappa_0$ is given by~\eqref{eq:kappa_0}.
\vspace{5mm}

\begin{lemma}\label{lemma:spectral_analysis}
    For the operator $A_{\kappa}$ defined by \eqref{eq:A_kappa}-\eqref{eq:domain_A_kappa}, the following statements hold:
    \begin{itemize}
    \item[(i)] The growth bound $\omega_0(\kappa)$ is continuous and increasing as a function of $\kappa\in\mathbb{R}$. 
    \item[(ii)] If $\kappa=u_0$, then
    \[
    \omega_0(\kappa)=-\frac{\pi^2}{4Pe}-\frac{P_e}{4}-r.
    \]
    \item[(iii)] If $\kappa<\bar{\kappa}$, $\kappa\neq u_0$, then
    \[
    \omega_0(\kappa)=-\frac{\mu^2}{Pe}-\frac{P_e}{4}-r,
    \]
    where $\mu$ is the minimal positive solution to the equation
    \begin{equation}\label{eq:negative_eigenvalues}
        \tan(\mu) = \frac{4\mu Pe(1-\kappa)}{Pe^2(2\kappa-1)+4\mu^2}.
    \end{equation}
    \item[(iv)] If $\kappa=\bar{\kappa}$, then
    \[
    \omega_0(\kappa)=-\frac{P_e}{4}-r.
    \]
    \item[(v)] If $\kappa=1$, then
    \[
     \omega_0(\kappa)=-r.
    \] 
    \item[(vi)] If $\kappa>\bar{\kappa}$, $\kappa\neq1$, then
    \[
    \omega_0(\kappa)=\frac{\mu^2}{Pe}-\frac{P_e}{4}-r,
    \]
    where $\mu$ is the maximal positive solution to the equation
    \begin{equation}\label{eq:positive_eigenvalues}
        e^{2\mu} = \frac{(2\mu+Pe(2\kappa-1))(2\mu-Pe)}{(2\mu-Pe(2\kappa-1))(2\mu+Pe)}.
    \end{equation}
    In particular, if $\kappa=\kappa_0>1$, then
    \[
     \omega_0(\kappa)=0.
    \]
    \end{itemize}
\end{lemma}
\begin{proof}
The boundedly invertible linear operator $P:\mathrm{L}^2(0,1)\to \mathrm{L}^2(0,1)$ defined by $(Pc)(\zeta) = e^{-\frac{Pe}{2}\zeta}c(\zeta),\, c\in \mathrm{L}^2(0,1),\,\zeta\in[0,1]$, induces the (positivity preserving) similarity transformation
\begin{equation*}
    PA_{\kappa}P^{-1} = \tilde{A}_{\kappa} := \frac{1}{Pe}\frac{{\rm  d}^2}{{\rm  d}\zeta^2} - \left(\frac{Pe}{4}+r\right)I
\end{equation*}
on the domain
\begin{equation*}
    \begin{aligned}
    D(\tilde{A}_{\kappa}) = \biggl\{x\in \mathrm{H}^2(0,1) : &\frac{{\rm  d}  x}{{\rm  d}  \zeta}(0)+Pe\left(\kappa-\frac{1}{2}\right)x(0) =  0, \\
    &\frac{{\rm  d}  x}{{\rm  d}  \zeta}(1)+\frac{Pe}{2}x(1) =0 \biggr\}.
    \end{aligned}
\end{equation*}
Therefore, the growth bound $\omega_0(\kappa)$ of the $C_0$-semigroup generated by $A_{\kappa}$ satisfies
\begin{equation}\label{eq:growth_bound_shift}
    \omega_0(\kappa) = \theta_0(\kappa)-\frac{Pe}{4}-r,
\end{equation}
where $\theta_0(\kappa)$ is the growth bound, or equivalently the maximal eigenvalue, of the shifted operator
\begin{equation}\label{eq:D_kappa}
    D_{\kappa} := \frac{1}{Pe}\frac{{\rm  d}^2}{{\rm  d}\zeta^2}, \quad D(D_{\kappa}) = D(\tilde{A}_{\kappa})
\end{equation}
(since $-D_{\kappa}$, with domain $D(-D_{\kappa}):=D(D_{\kappa})$, is a Sturm-Liouville operator). In the following, it suffices then to study the eigenvalues of $D_{\kappa}$, and more particularly $\theta_0(\kappa)$. By invoking arguments similar to those used in the proof of \cite[Lemma 3.3]{Yevgenieva2026} where $\cmcal{L}=-D_{\kappa}$, we obtain that 
\begin{itemize}
    \item[\textbf{(A)}] $\theta_0(\kappa)$ is a continuous and increasing function of $\kappa\in\mathbb{R}$,
    \item[\textbf{(B)}] $0$ is an eigenvalue of $D_{\kappa}$ if, and only if, $\kappa=\bar{\kappa}$,
    \item[\textbf{(C)}] if $\kappa>\bar{\kappa}$, then $D_{\kappa}$ has at least one positive eigenvalue,
    \item[\textbf{(D)}] if $\kappa\leq1/2$, then $D_{\kappa}$ has no positive eigenvalue,
    \item[\textbf{(E)}] if $\theta$ is a positive eigenvalue of $D_{\kappa}$, then $\theta=Pe/4$ and $\kappa=1$, or $\theta=\mu^2/Pe$ where $\mu$ is a positive root of \eqref{eq:positive_eigenvalues} and $\kappa\neq1$,
    \item[\textbf{(F)}] if $\theta$ is a negative eigenvalue of $D_{\kappa}$, then $\theta=-\mu^2/Pe$ where $\mu$ is a positive root of \eqref{eq:negative_eigenvalues} or there exists $n\in\mathbb{N}$ such that $\kappa=u_n$ and $\mu=\left(n+\frac{1}{2}\right)\pi$.
\end{itemize}
The statement (i) follows directly from \textbf{(A)} and \eqref{eq:growth_bound_shift}. The property  \textbf{(C)} implies that $\theta_0(\kappa)>0$ for $\kappa>\bar{\kappa}$; together with \textbf{(E)} where $\theta=\theta_0(\kappa)$ and \eqref{eq:growth_bound_shift}, this proves (v) and (vi).\\ 
Suppose now that $\theta_0(\bar{\kappa})\neq0$, that is $\theta_0(\bar{\kappa})>0$, in view of \textbf{(B)}. Observe that \textbf{(B)} and \textbf{(D)} lead to the inequality $\theta_0(\kappa)<0$ for $\kappa\leq1/2$. Hence, because of \textbf{(A)}, there would exist $\tilde{\kappa}\in(1/2,\bar{\kappa})$ such that $\theta_0(\tilde{\kappa})=0$, which leads to a contradiction with \textbf{(B)}. Consequently, $\theta_0(\bar{\kappa})=0$ and (iv) follows from \eqref{eq:growth_bound_shift}.\\ 
By invoking \textbf{(A)} again, observe that $\theta_0(\kappa)<\theta_0(\bar{\kappa})=0$ when $\kappa<\bar{\kappa}$. In addition, for all $n\in\mathbb{N}$, it can be shown that $\left(n+\frac{1}{2}\right)\pi$ is smaller than the minimal positive root of \eqref{eq:negative_eigenvalues} with $\kappa=u_n$ if, and only if, $n=0$. If $\kappa<\bar{\kappa}$, the latter result, when combined with \eqref{eq:growth_bound_shift} and \textbf{(F)} for $\theta=\theta_0(\kappa)<0$, proves (ii) and (iii).\\
\end{proof}

\begin{figure}[thpb]
      \centering
      \vspace{-0.4cm}
      \includegraphics[scale=0.55]{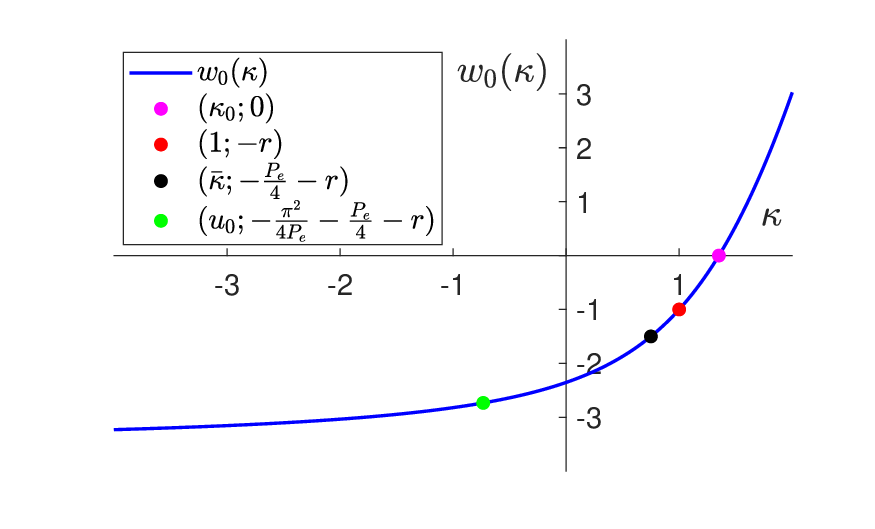}
      \vspace{-0.7cm}
      \caption{Growth bound $\omega_0(\kappa)$ as a function of $\kappa$ for $Pe=2$ and $r=1$}
      \label{fig:growth_bound}
   \end{figure}

Figure \ref{fig:growth_bound} represents $\omega_0(\kappa)$, explicitly given in Lemma \ref{lemma:spectral_analysis} as a function of $\kappa$, for $Pe=2$ and $r=1$.
   
Now we are able to prove the following result that provides an analytical expression of the solution of $\mathbf{P_{\alpha}}$.\\

\begin{theorem}\label{theorem:solution_pbm_opti}
    Fix~$\alpha>0$. The $\mathbf{P_{\alpha}}$ problem~\eqref{eq:pbm_opti} for the reaction-convection-diffusion system \eqref{eq:system}-\eqref{eq:boundary_cdts} with a state feedback of the form~\eqref{eq:feedback_exact} is feasible, i.e. there exists $\kappa\in\mathbb{R}$ such that \eqref{eq:pbm_opti_stability}, \eqref{eq:pbm_opti_system_positivity} and~\eqref{eq:pbm_opti_input_positivity} hold. Moreover, $\mathbf{P_{\alpha}}$ has a unique solution $\kappa^*$ given by
    \begin{equation*}
        \kappa^* = \left\{
        \begin{array}{ll}
            \kappa_{\alpha}, &\text{ if } \alpha\leq-\omega_0(0),\\
            0, &\text{ otherwise, }
        \end{array}\right.
    \end{equation*}
    where, for all $\alpha\in (0,-\omega_0(0)]$, $\kappa_{\alpha}$ is the feedback gain such that the growth constant of the resulting closed-loop semigroup $\omega_0(\kappa_{\alpha}) = -\alpha$.\\
\end{theorem}
\begin{proof}
    Let $\alpha>0$. Using arguments similar to the one used in \cite[Lemma 5.1]{Laabissi2001}, we show that $(T_{D_{\kappa}}(t))_{t\geq0}$, and so $(T_{A_{\kappa}}(t))_{t\geq0}$, is positive for all $\kappa\in\mathbb{R}$. Invoking Lemma \ref{lemma:spectral_analysis}, and in particular the fact that $\omega_0(0)<0$, $\omega_0(\kappa_0)=0$ with $\kappa_0>0$ and the function $\omega_0(\kappa)$ is continuous and increasing  on $[0,\kappa_0]$, shows that the combination of $\kappa\geq0$ and $-\alpha\leq\omega_0(\kappa)<0$ is equivalent to 
    \begin{equation*}
        \left\{
        \begin{array}{ll}
            \kappa_{\alpha}\leq\kappa<\kappa_0, &\text{ if } \alpha\leq-\omega_0(0),\\
            0\leq\kappa<\kappa_0, &\text{ otherwise. }
        \end{array}\right.
    \end{equation*}
    where, for all $\alpha\leq-\omega_0(0)$,  $\kappa_{\alpha}$ $(\geq0)$ is the unique solution of $\omega_0(\kappa_{\alpha}) = -\alpha$. 
    Therefore, the $\mathbf{P_{\alpha}}$ problem is feasible and its feasible set is given by $[\kappa_{\alpha},\kappa_0)$ when $\alpha\leq-\omega_0(0)$, and by $[0,\kappa_0)$ otherwise. In either case, the feasible set is contained in $\mathbb{R}_+$. Since the function $\kappa\mapsto\kappa^2$ is increasing on $\mathbb{R}_+$, its minimum over the feasible set is attained at the lower endpoint of the feasible interval. This concludes the proof.\\
\end{proof}

\begin{remark}
    The precise value of $\omega_0(0)$ can be obtained straightforwardly thanks to Lemma \ref{lemma:spectral_analysis}, as well as the one of $\kappa_{\alpha}$ for any $\alpha\in(0,-\omega_0(0)]$. More precisely, 
    \begin{itemize}
        \item if $\alpha\in\left(\frac{Pe}{4}+r,-\omega_0(0)\right]$, $\alpha\neq\frac{\pi^2}{4Pe}+\frac{P_e}{4}+r$, then
        \begin{equation*}
            \left\{
            \begin{aligned}
                \kappa_{\alpha} &= \frac{4\mu Pe+\tan(\mu)(Pe^2-4\mu^2)}{4\mu Pe+2\tan(\mu)Pe^2},\\
                \mu&=\sqrt{\left(\alpha-\frac{Pe}{4}-r\right)Pe},
            \end{aligned}\right.
        \end{equation*}
        \item if $\alpha=\frac{\pi^2}{4Pe}+\frac{P_e}{4}+r\leq-\omega_0(0)$, then $\kappa_{\alpha}=u_0$,
        \item if $\alpha=\frac{P_e}{4}+r$, then $\kappa_{\alpha}=\bar{\kappa}$,
        \item if $\alpha\in\left[0,\frac{Pe}{4}+r\right)$, $\alpha\neq r$, then
        \begin{equation*}
            \left\{
            \begin{aligned}
                \kappa_{\alpha} &= \frac{4\mu Pe+\tanh(\mu)(Pe^2+4\mu^2)}{4\mu Pe+2\tanh(\mu)Pe^2},\\
                \mu&=\sqrt{\left(-\alpha+\frac{Pe}{4}+r\right)Pe},
            \end{aligned}\right.
        \end{equation*}
        \item if $\alpha=r$, then $\kappa_{\alpha}=1$.
    \end{itemize}
    By monotonicity, observe that the condition $\frac{\pi^2}{4Pe}+\frac{P_e}{4}+r\leq-\omega_0(0)$ is equivalent to $0\leq\kappa_{\frac{\pi^2}{4Pe}+\frac{P_e}{4}+r}=u_0$, i.e. $Pe\geq\pi$. 
\end{remark}

\section{OBSERVER-BASED FEEDBACK}\label{section:observer_based_feedback}
In section \ref{section:exact_feedback}, exact state feedback laws of the form $c_{in}(\tau)=\kappa c(0,\tau)$ were developed to regulate the state dynamics of system \eqref{eq:system}-\eqref{eq:boundary_cdts} in accordance with desired process requirements. The practical implementation of these feedback laws, however, relies on the availability of sufficiently accurate knowledge of the state $c(\cdot,\tau)$ at every time $\tau\geq0$. In the present setting, only the concentration of the reactant at the outlet is measurable; see \eqref{eq:output}. Consequently, to overcome this limitation, we replace the exact state feedback law \eqref{eq:feedback_exact} by an observer-based feedback strategy for the diffusion system \eqref{eq:system}-\eqref{eq:boundary_cdts}, defined as follows:
\begin{equation}\label{eq:observer_feedback}
    c_{in}(\tau) = \kappa \tilde{c}(0,\tau)\,, \quad \tau\geq 0\,,
\end{equation}
where $\kappa \in \mathbb{R}$ still denotes the feedback gain and $\tilde{c}(\cdot,t)$ is the estimated state generated by the observer dynamics
\begin{equation}\label{eq:observer}
    \frac{\partial \tilde{c}}{\partial \tau}(\zeta,\tau) = \frac{1}{Pe}\frac{\partial^2 \tilde{c}}{\partial {\zeta}^2}(\zeta,\tau)-\frac{\partial \tilde{c}}{\partial {\zeta}}(\zeta,\tau)-r\tilde{c}(\zeta,\tau),
\end{equation}
with boundary conditions 
\begin{equation}\label{eq:observer_boundary_cdts}
\left\{
\begin{array}{lcl}
    \dfrac{1}{Pe}\dfrac{\partial \tilde{c}}{\partial \zeta}(0,\tau)-\tilde{c}(0,\tau) & = & -c_{in}(\tau),\\
    \dfrac{1}{Pe}\dfrac{\partial \tilde{c}}{\partial \zeta}(1,\tau) & = & \gamma(\tilde{y}(\tau)-y(\tau))\,,
\end{array}
\right.
\end{equation}
where the measured output $y$ is given by \eqref{eq:output}, and where
\begin{equation}\label{eq:estimated_output}
    \tilde{y}(\tau):=\tilde{c}(1,\tau),\quad t\geq 0,
\end{equation}
denotes the estimated output, and $\gamma \in \mathbb{R}$ is the observer gain, to be designed. The dynamics of the state estimation error $e:=\tilde{x}-x$ is described by the following PDE:
\begin{equation*}
    \frac{\partial e}{\partial \tau}(\zeta,\tau) = \frac{1}{Pe}\frac{\partial^2 e}{\partial {\zeta}^2}(\zeta,\tau)-\frac{\partial e}{\partial {\zeta}}(\zeta,\tau)-re(\zeta,\tau),
\end{equation*}
with Robin’s boundary conditions 
\begin{equation*}
\left\{
\begin{array}{lcl}
    \dfrac{1}{Pe}\dfrac{\partial e}{\partial \zeta}(0,\tau)-e(0,\tau) & = & 0,\\
    \dfrac{1}{Pe}\dfrac{\partial e}{\partial \zeta}(1,\tau) -\gamma e(1,\tau)& = & 0,
\end{array}
\right.
\end{equation*}
Through this whole section, we will make the following \textbf{Standing Assumption:} 
$\kappa \; \neq \; \gamma$.\\

As explained in section \ref{section:model_problem}, one objective of any designed control strategy under consideration is to preserve the exponential convergence of the reactant trajectories to zero. Moreover, when the reactant dynamics is combined with a state observer, the resulting state estimation error should satisfy the same property. 
It turns out that both convergence properties are consequences of the exponential stability of the dynamics of $c$ and $e$ together, i.e. of the $(c,e)$-dynamics. Lemma \ref{lemma:stability} provides necessary and sufficient conditions on the feedback and observer gains for the $(c,e)$-dynamics to be exponentially stable. \\

\begin{lemma}\label{lemma:stability}
    The $(c,e)$-dynamics, i.e. 
    \begin{equation*}
        \begin{pmatrix}
            \dot{c}(\tau)\\\dot{e}(\tau)
        \end{pmatrix}
        =
        \mathcal{A}_{\kappa\gamma}\begin{pmatrix}
            c(\tau)\\e(\tau)
        \end{pmatrix}, 
    \end{equation*}
    where the operator $\mathcal{A}_{\kappa\gamma}:D(\mathcal{A}_{\kappa\gamma})\subset\mathrm{L}^2(0,1)\times\mathrm{L}^2(0,1)\rightarrow \mathrm{L}^2(0,1)\times\mathrm{L}^2(0,1)$ is defined by
    \begin{equation*}
        \mathcal{A}_{\kappa\gamma}=\begin{pmatrix}
            \frac{1}{Pe}\frac{{\rm  d}^2}{{\rm  d}  \zeta^2}-\frac{{\rm  d}}{{\rm  d}  \zeta}-rI&0\\
            0&\frac{1}{Pe}\frac{{\rm  d}^2}{{\rm  d}  \zeta^2}-\frac{{\rm  d}}{{\rm  d}  \zeta}-rI
        \end{pmatrix}
    \end{equation*}
    on its domain
    \begin{equation*}
        \begin{aligned}
            D(\mathcal{A}_{\kappa\gamma}) = \biggl\{& X:=(c\hspace{0.2cm}e)^{\top}\in \mathrm{H}^2(0,1) \times \mathrm{H}^2(0,1):\\ &\frac{1}{Pe}\dfrac{{\rm  d}}{{\rm  d}  \zeta}\begin{pmatrix}c\\e\end{pmatrix}(0) = \begin{pmatrix}
            1-\kappa&-\kappa\\0&1\end{pmatrix}\begin{pmatrix}c\\e\end{pmatrix}(0),\\ &\frac{1}{Pe}\dfrac{{\rm  d}}{{\rm  d}  \zeta}\begin{pmatrix}c\\e\end{pmatrix}(1) = \begin{pmatrix}
            0&0\\0&\gamma\end{pmatrix}\begin{pmatrix}c\\e\end{pmatrix}(1) \biggl\}
        \end{aligned}
    \end{equation*}
    is exponentially stable if, and only if, $\kappa<\kappa_0$ and $\gamma<\kappa_0$.\\
\end{lemma}
\begin{proof}
    Let $\kappa,\gamma \in\mathbb{R}$. Let us define the operator $$\mathcal{D}_{\kappa\gamma} := \mathcal{P}\mathcal{A}_{\kappa\gamma}\mathcal{P}^{-1}+\left(\frac{Pe}{4}+r\right)I$$ with domain $D(\mathcal{D}_{\kappa\gamma}) := D(\mathcal{P}\mathcal{A}_{\kappa\gamma}\mathcal{P}^{-1})$, where $\mathcal{P}$ is the similarity transformation defined by $\mathcal{P}X:=(Pc\hspace{0.2cm}Pe)^{\top}$ for all $X=(c\hspace{0.2cm}e)^{\top}\in \mathrm{L}^2(0,1) \times \mathrm{L}^2(0,1)$ with $P$ introduced in the proof of Lemma \ref{lemma:spectral_analysis}. Observe that $\mathcal{D}_{\kappa\gamma}$ is a 2D diffusion operator with a domain $D(\mathcal{D}_{\kappa\gamma})$ given by 
        \begin{equation*}
        \begin{aligned}
            D(\mathcal{D}_{\kappa\gamma}) = \biggl\{& X:=(c\hspace{0.2cm}e)^{\top}\in \mathrm{H}^2(0,1) \times \mathrm{H}^2(0,1):\\ &\frac{1}{Pe}\dfrac{{\rm  d}}{{\rm  d}  \zeta}\begin{pmatrix}c\\e\end{pmatrix}(0) = \begin{pmatrix}
            \frac{1}{2}-\kappa&-\kappa\\0&\frac{1}{2}\end{pmatrix}\begin{pmatrix}c\\e\end{pmatrix}(0),\\ &\frac{1}{Pe}\dfrac{{\rm  d}}{{\rm  d}  \zeta}\begin{pmatrix}c\\e\end{pmatrix}(1) = \begin{pmatrix}
            -\frac{1}{2}&0\\0&\gamma-\frac{1}{2}\end{pmatrix}\begin{pmatrix}c\\e\end{pmatrix}(1) \biggl\},
        \end{aligned}
    \end{equation*}
    i.e. $D(\mathcal{D}_{\kappa\gamma})$ is defined by 2D Robin's boundary conditions that are described by upper triangular matrices. Consider now the eigenvalue problem 
    \begin{equation}\label{eq:eigenvalue_problem}
        \mathcal{D}_{\kappa\gamma}\Psi = \theta\Psi,\quad\theta\in\mathbb{C},\:\Psi\in D(\mathcal{D}_{\kappa\gamma}).
    \end{equation}
    By taking advantage of the triangular structure of the boundary conditions in $D(\mathcal{D}_{\kappa\gamma})$, a standard analysis of \eqref{eq:eigenvalue_problem} shows that $\sigma_p(\mathcal{D}_{\kappa\gamma}) = \sigma_p(D_{\kappa})\cup\sigma_p(E_{\gamma})$, where $D_{\kappa}$ is defined by \eqref{eq:D_kappa} and $E_{\gamma}:=\frac{1}{Pe}\frac{{\rm  d}^2}{{\rm  d}\zeta^2}$ with domain
    \begin{equation*}
        \begin{aligned}
            D(E_{\gamma}) := \biggl\{e\in \mathrm{H}^2(0,1) : &\frac{{\rm  d}  e}{{\rm  d}  \zeta}(0)-\frac{Pe}{2}e(0) =  0, \\
            &\frac{{\rm  d}  e}{{\rm  d}  \zeta}(1)-Pe\left(\gamma-\frac{1}{2}\right)e(1) =0 \biggr\}.
        \end{aligned}
    \end{equation*}
    Now, by introducing the invertible map $S:\mathrm{L}^2(0,1)\to \mathrm{L}^2(0,1)$ defined by $(Se)(\zeta) = e(1-\zeta),\, e\in \mathrm{L}^2(0,1),\,\zeta\in[0,1],$
    we find that $SE_{\gamma}S^{-1} = D_{\gamma}$ with $D(SE_{\gamma}S^{-1}) = D(D_{\gamma})$, and so $\sigma_p(E_{\gamma})=\sigma_p(D_{\gamma})$. Therefore, 
    $\sigma_p(\mathcal{D}_{\kappa\gamma}) = \sigma_p(D_{\kappa})\cup\sigma_p(D_{\gamma})$. Moreover, the standing assumption implies that $\sigma_p(D_{\kappa})\cap\sigma_p(D_{\gamma})\neq\emptyset$. From the spectral analysis performed in Lemma~\ref{lemma:spectral_analysis}, we then get that $\mathcal{D}_{\kappa\gamma}$ is a Riesz spectral operator~\cite[Definition 3.2.6]{Curtain_Zwart2020} with simple eigenvalues and a finite spectral bound. It follows, by \cite[Theorem 3.2.8]{Curtain_Zwart2020}, that $\mathcal{D}_{\kappa\gamma}$ generates a $C_0$-semigroup $(T_{\mathcal{D}_{\kappa\gamma}}(t))_{t\geq0}$ with growth bound  $\omega_0(T_{\mathcal{D}_{\kappa\gamma}}) = s(\mathcal{D}_{\kappa\gamma})$. Thus, $\mathcal{A}_{\kappa\gamma}$ generates a $C_0$-semigroup $(T_{\mathcal{A}_{\kappa\gamma}}(t))_{t\geq0}$ such that $T_{\mathcal{A}_{\kappa\gamma}}(t) = \mathcal{P}^{-1}e^{-(\frac{Pe}{4}+r)t}T_{\mathcal{D}_{\kappa\gamma}}(t)\mathcal{P}$ and
        \begin{align}
            \omega_0(T_{\mathcal{A}_{\kappa\gamma}})& =\omega_0(T_{\mathcal{D}_{\kappa\gamma}})-\frac{Pe}{4}-r\notag\\
            &= s(\mathcal{D}_{\kappa\gamma})-\frac{Pe}{4}-r\notag\\
            &= \max(\theta_0(\kappa),\theta_0(\gamma))-\frac{Pe}{4}-r\notag\\
            &=\max(\omega_0(\kappa),\omega_0(\gamma)),\label{eq:growth_bound_c_e}
        \end{align}
    where, for all $\delta\in\mathbb{R}$, $\theta_0(\delta)$ and $\omega_0(\delta)$ are defined in the proof of Lemma~\ref{lemma:spectral_analysis}. Invoking Lemma~\ref{lemma:spectral_analysis} with \eqref{eq:growth_bound_c_e}, allows to conclude that $(T_{\mathcal{A}_{\kappa\gamma}}(t))_{t\geq0}$ is exponentially stable if, and only if, $\kappa<\kappa_0$ and $\gamma<\kappa_0$. 
\end{proof}

\vspace{5mm}

In the proof of Lemma \ref{lemma:stability}, it is actually shown that the error dynamics is equivalent to the closed-loop dynamics \eqref{eq:closed_loop}-\eqref{eq:domain_A_kappa} studied in section \ref{section:exact_feedback}, with $\kappa$ replaced by $\gamma$. We also establish the validity of the so-called separation principle in our setting.

Let us now turn our attention to the positivity property. When regulating any positive system with an observer-based state feedback, it is necessary to design feedback and observers gains that produce states maintaining this intrinsic property. Positivity of estimated states is also highly recommended; otherwise that would have no physical meaning. The positivity of the $(c,\tilde{c})$-dynamics can be characterized in terms of the positivity of a simpler system, as established by the following result. \\

\begin{lemma}\label{lemma:positivity}
    The $(c,\tilde{c})$-dynamics (i.e. \eqref{eq:system}-\eqref{eq:output} combined with \eqref{eq:observer}-\eqref{eq:estimated_output} and \eqref{eq:observer_feedback}) is positive if, and only if, the auxiliary system 
        \begin{equation}\label{eq:systemDiff_c_ctilde}
            \begin{pmatrix}
                \frac{\partial c}{\partial \tau}(\zeta,\tau) \\
                \frac{\partial \tilde{c}}{\partial \tau}(\zeta,\tau)
            \end{pmatrix}
            =
            \frac{1}{Pe}
            \begin{pmatrix}
                \frac{\partial^2 c}{\partial {\zeta}^2}(\zeta,\tau)\\
                \frac{\partial^2 \tilde{c}}{\partial {\zeta}^2}(\zeta,\tau)
            \end{pmatrix}
        \end{equation}
        with boundary conditions 
        \begin{equation}\label{eq:boundary_cdts_Diff_c_ctilde}
        \left\{
            \begin{aligned}
                &\frac{1}{Pe}\begin{pmatrix}
                    \frac{\partial c}{\partial {\zeta}}(0,\tau)\\
                    \frac{\partial \tilde{c}}{\partial {\zeta}}(0,\tau)
                \end{pmatrix}
                =
                \begin{pmatrix}
                    \frac{1}{2}&-\kappa\\0&\frac{1}{2}-\kappa
                \end{pmatrix}
                \begin{pmatrix}
                    c(0,\tau)\\
                    \tilde{c}(0,\tau)
                \end{pmatrix},
                \\
                &\frac{1}{Pe}\begin{pmatrix}
                    \frac{\partial c}{\partial {\zeta}}(1,\tau)\\
                    \frac{\partial \tilde{c}}{\partial {\zeta}}(1,\tau)
                \end{pmatrix}
                =
                \begin{pmatrix}
                    -\frac{1}{2}&0\\-\gamma&\gamma-\frac{1}{2}
                \end{pmatrix}
                \begin{pmatrix}
                    c(1,\tau)\\
                    \tilde{c}(1,\tau)
                \end{pmatrix},
            \end{aligned}
        \right.
        \end{equation}
        is positive.\\
\end{lemma}

\begin{proof}
    First, note that the operator associated with the $(c,\tilde{c})$-dynamics, denoted by $\mathcal{\tilde{A}}_{\kappa\gamma}$, generates a $C_0$-semigroup because so does $\mathcal{A}_{\kappa\gamma}$. Then, it suffices to observe that $$\mathcal{P}\mathcal{\tilde{A}}_{\kappa\gamma}\mathcal{P}^{-1}+\left(\frac{Pe}{4}+r\right)I=\mathcal{\tilde{D}}_{\kappa\gamma}, $$ and $$D(\mathcal{P}\mathcal{\tilde{A}}_{\kappa\gamma}\mathcal{P}^{-1})=D(\mathcal{\tilde{D}}_{\kappa\gamma}), $$ where the positivity preserving operator $\mathcal{P}$ is defined in the proof of Lemma \ref{lemma:stability} and $\mathcal{\tilde{D}}_{\kappa\gamma}$ with domain  $D(\mathcal{\tilde{D}}_{\kappa\gamma})$ is the operator associated with the dynamics \eqref{eq:systemDiff_c_ctilde}-\eqref{eq:boundary_cdts_Diff_c_ctilde}.
\end{proof}

\vspace{5mm}

The following result provides necessary and sufficient conditions on the feedback and observer gains, $\kappa$ and $\gamma$, respectively, ensuring that the $(c,e)$-dynamics is exponentially stable and that the $(c,\tilde{c})$-dynamics is positive.\\
\begin{theorem}\label{theorem:stability+positivity}
    The $(c,e)$-dynamics is exponentially stable and the $(c,\tilde{c})$-dynamics is positive if, and only if, $\kappa<\kappa_0$, $\gamma<\kappa_0$, and the auxiliary system \eqref{eq:systemDiff_c_ctilde}-\eqref{eq:boundary_cdts_Diff_c_ctilde} is positive.\\
    
\end{theorem}
\begin{proof}
    The result follows directly by combining Lemma~\ref{lemma:stability} and Lemma~\ref{lemma:positivity}. 
\end{proof}

\vspace{5mm}

For an uncontrolled infinite-dimensional system that is not exponentially stable, results in the literature establish that it is impossible to positively stabilize it with an observer-based feedback resulting from a positive observer; see \cite[Theorem 3.1]{Binid2021}, \cite[Theorem 5.1]{Rami2011} or \cite[Subsection 5.1]{PiengeonWinkin2026}. However, it is worth recalling that the system \eqref{eq:system}-\eqref{eq:boundary_cdts} with $c_{in}=0$ is already exponentially stable. Consequently, the positivity and stability requirements may potentially be satisfied together, i.e. values of  $\kappa<\kappa_0$ and $\gamma<\kappa_0$ can be found such that the auxiliary system \eqref{eq:systemDiff_c_ctilde}-\eqref{eq:boundary_cdts_Diff_c_ctilde} is positive. This is strengthened in the following section with Proposition \ref{proposition:positivity_discretized_system} and numerical simulations. 

\section{NUMERICAL SIMULATIONS}\label{section:numerical_simulations}
In this section, the goal is to get an educated guess for the values of $\kappa$ and $\gamma$ that would ensure that the auxiliary system \eqref{eq:systemDiff_c_ctilde}-\eqref{eq:boundary_cdts_Diff_c_ctilde} be positive. To this end, let us discretize this system. Consider $n\in\mathbb{N}$ ($n\neq0,1$) discretization points $\zeta_i$, $i=1,\dots,n$, with $\zeta_1=0$ and $\zeta_n=1$. The discretization step, denoted $\Delta \zeta$, is given by~$\Delta \zeta = 1/(n-1)$. Using second-order central differences about $\zeta_i$ for $\frac{\partial^2 c}{\partial {\zeta}^2}(\zeta_i,\tau)$ and $\frac{\partial^2 \tilde{c}}{\partial {\zeta}^2}(\zeta_i,\tau)$, $i=1,\dots,n$, first-order backward differences about $0$ for $\frac{\partial c}{\partial {\zeta}}(0,\tau)$ and $\frac{\partial \tilde{c}}{\partial {\zeta}}(0,\tau)$ and first-order forward differences about $1$ for $\frac{\partial c}{\partial {\zeta}}(1,\tau)$ and $\frac{\partial \tilde{c}}{\partial {\zeta}}(1,\tau)$, we obtain a $2n$-dimensional discretization of \eqref{eq:systemDiff_c_ctilde}-\eqref{eq:boundary_cdts_Diff_c_ctilde} given by \begin{equation}\label{eq:discretized_system}
    \dot{X}^{(2n)}(\tau) = \mathcal{\tilde{D}}^{(2n)}_{\kappa\gamma}X^{(2n)}(\tau), 
\end{equation}
where $X^{(2n)}(\cdot)\approx [c(\zeta_1,\cdot) \cdots c(\zeta_n,\cdot) \: \tilde{c}(\zeta_1,\cdot) \cdots \tilde{c}(\zeta_n,\cdot)]^{\top}$ and $\mathcal{\tilde{D}}^{(2n)}_{\kappa\gamma}=:[d_{ij}]\in\mathbb{R}^{2n\times 2n}$ is a matrix whose nonzero off-diagonal entries are $d_{i+1,i}=d_{i,i+1}=1/(\Delta\zeta^2 Pe)$, $i=1,\dots,2n-1$, $i\neq n$, $d_{1,n+1}=\kappa/\Delta\zeta$ and $d_{2n,n}=-\gamma/\Delta\zeta$. \\
\begin{proposition}\label{proposition:positivity_discretized_system}
    For any $n\in\mathbb{N}$ ($n\neq0,1$), the spatial discretization~\eqref{eq:discretized_system} of the auxiliary system \eqref{eq:systemDiff_c_ctilde}-\eqref{eq:boundary_cdts_Diff_c_ctilde} is positive if, and only if, $\kappa\geq0$ and $\gamma\leq0$. \\
\end{proposition}
\begin{proof}
     It suffices to observe that $\mathcal{\tilde{D}}^{(2n)}_{\kappa\gamma}$ is Metzler, or equivalently \cite[Part I, Theorem 2]{Farina2000}, the finite-dimensional system \eqref{eq:discretized_system} is positive if, and only if, $\kappa\geq0$ and $\gamma\leq0$. \\
\end{proof}

\begin{figure}[thpb]
      \centering
      \vspace{-0.4cm}
      \includegraphics[scale=0.5]{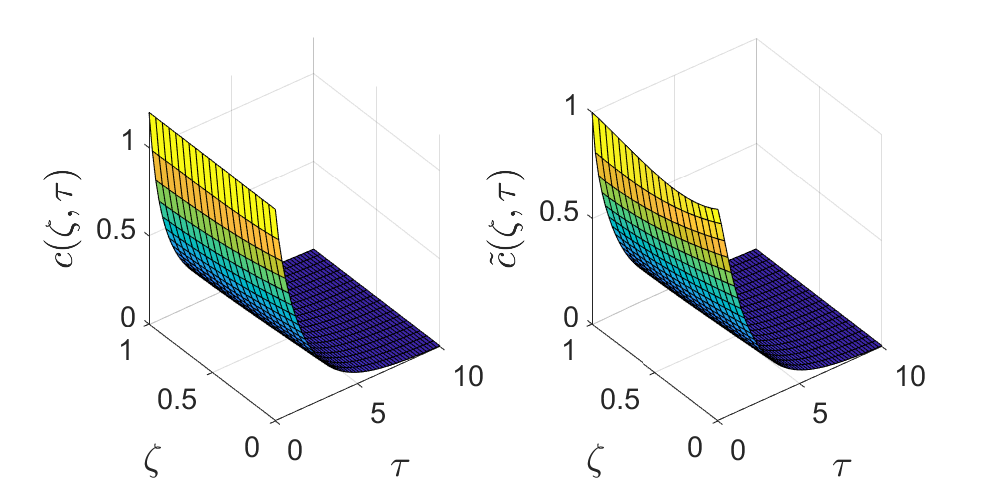}
      \vspace{-0.2cm}
      \caption{State and estimated state trajectories for $Pe=2$, $r=1$, $\kappa=1.2<\kappa_0=1.3516$, and $\gamma=-0.2$}
      \label{fig:trajectories}
   \end{figure}

The combination of Theorem \ref{theorem:stability+positivity} and Proposition \ref{proposition:positivity_discretized_system} with the positivity of $\kappa_0$ suggests that taking $0\leq\kappa<\kappa_0$ and $\gamma\leq0$ are necessary and sufficient conditions for the $(c,e)$-dynamics to be exponentially stable and the $(c,\tilde{c})$-dynamics to be positive. If this equivalence is verified, then such gains yield a nonnegative inlet concentration $c_{in}$, in accordance with physical requirements. Figure \ref{fig:trajectories} shows the state and estimated state trajectories of the $(c,\tilde{c})$-dynamics obtained from the nonnegative initial conditions $c_0(\zeta)=\kappa$ and $\tilde{c}_0(\zeta)=e^{Pe(1-\kappa)\zeta(\zeta-1)^2+Pe\gamma(1-\kappa)\zeta^2(\zeta-1)}$, $\zeta\in[0,1]$, with the parameters $Pe=2, r=1, \kappa=1.2<\kappa_0=1.3516$, and $\gamma=-0.2$. Both trajectories converge to zero while staying nonnegative. 

\section{CONCLUDING REMARKS}\label{section:conclusion}
In this paper, we have explained how to control the reactant concentration in an isothermal axial dispersion tubular reactor according to relevant process requirements. Three of them were maintained in the whole analysis: convergence to zero and nonnegativity of the state trajectories, and nonnegativity of the inlet reactant concentration. First, a closed-form solution for an optimal feedback design problem was derived. This solution may be viewed as a trade-off between a sufficiently slow decay rate and sufficiently low feedback static gains, while satisfying the three requirements. Then, we have explained how to obtain an educated guess for the positive stabilization of an observer-based closed-loop system.\\ 
An ongoing research consists in the validation of this educated guess with, for example, the study of the resolvent operator associated with the operator $\mathcal{\tilde{A}}_{\kappa\gamma}$. As possible perspectives, it could be relevant to analyze the performance of the designed control on the product concentration of the considered model. Finally, we mention the positive control of nonlinear tubular reactors, as an outlook for future research.






\bibliographystyle{IEEEtran}
\bibliography{IEEEabrv,mybibfile}

\end{document}